\documentclass[10pt]{amsart}

\title[Remarks on the thin obstacle problem and constrained Ginibre ensembles]{Remarks on the thin obstacle problem and constrained Ginibre ensembles}
\author{Aram L. Karakhanyan}
\address{School of Mathematics, The University of Edinburgh, Peter Tait Guthrie Road, EH9 3FD, Edinburgh, UK}
\email{ aram6k@gmail.com}

\usepackage{amsmath,amsfonts,amssymb,amsthm,upref}

\usepackage{subfigure}

\usepackage{hyperref}

\usepackage{enumitem}

\usepackage{verbatim}

\usepackage{skak} 

\usepackage{mathrsfs}
\usepackage{esint}

\usepackage{graphicx, geometry}
\hypersetup{
    bookmarks=true,         
    unicode=true,          
    pdftoolbar=true,        
    pdfmenubar=true,        
    pdffitwindow=false,     
    pdfstartview={FitH},    
    pdftitle={My title},    
    pdfauthor={Author},     
    pdfsubject={Subject},   
    pdfcreator={Creator},   
    pdfproducer={Producer}, 
    pdfkeywords={keyword1} {key2} {key3}, 
    pdfnewwindow=true,      
    colorlinks=true,       
    linkcolor=blue,          
    citecolor=blue,        
    filecolor=magenta,      
    urlcolor=cyan           
}

\renewenvironment{proof}[1][\proofname ]{{\noindent \bfseries #1. }}{\qed \bigskip }

\newcommand{\R}{{\mathbb R}}

\newcommand\dist{\operatorname{dist}} 

\newcommand{\e}{\varepsilon}
\newcommand{\supp}{\operatorname{supp}}

\renewcommand\H{\mathcal H}

\newcommand\p{\partial}

\def\Om{\Omega}
\def\na{\nabla}

\newcommand\M{\mathcal M}

\newcommand{\be}{\begin{equation}}
\newcommand{\ee}{\end{equation}}

\def\ua{U^{\mu_a}}
\newcommand\F[1]{\widehat{#1}}

\newtheorem{theorem}{Theorem}[section]

\newtheorem{cor}[theorem]{Corollary}

\newtheorem{lem}[theorem]{Lemma}

\newtheorem{remark}[theorem]{Remark}

\newtheoremstyle{named}{}{}{\itshape}{}{\bfseries}{.}{.5em}{\thmnote{#3 }#1}
\theoremstyle{named}

\numberwithin{equation}{section}


\begin{document}
\maketitle

\baselineskip=13pt    

\begin{abstract}
We consider the problem of constrained Ginibre ensemble  with prescribed portion of 
eigenvalues on a given curve $\Gamma\subset \R^2$ and relate it to a thin obstacle problem.
The key step in the proof is the $H^1$ estimate for the logarithmic potential of the equilibrium measure. 
The coincidence set has two components: one in $\Gamma$ and another one in $\R^2\setminus \Gamma$ which are well separated. 
Our main result here asserts that this obstacle problem is well posed in $H^1(\R^2)$ which 
improves previous results in $H^1_{loc}(\R^2)$.
\end{abstract}

\section{Introduction}
Let $\Gamma$ be a regular curve in $\R^2$ with locally finite length 
and $\mathcal M_a$  the set of all probability measures such that 
\be\label{def-M}
\mu(\Gamma)\ge a, 
\quad a\in (0,1).
\ee
 By an abuse of notation we let 
$\Gamma: \R\to \R^2$ be the arc-length parametrization of the curve such that 
\[|\dot{\Gamma} (t)|=1,\quad  t\in \R.\]
 In this paper we consider the minimizers of the energy 
\be
I[\mu]=\int\int\log\frac1{|x-y|}d\mu(x)d\mu(y)+\int Qd\mu
\ee
where  
$Q(x)$ is a given function such that  the weight function $w=e^{-Q}$ on $\R^2$ is admissible 
(see Definition 1.1 p.26 \cite{Saff}). This means that 
$w$ satisfies the following three conditions:
\begin{itemize}
\item[\bf(H1)] $w$ is upper semi-continuous;
\item[\bf(H2)] $\{w\in \R^2 \ s.t. \ w(z) > 0\}$ has positive capacity;
\item[\bf(H3)]  $|z|w(z) \to 0 $ as $|z| \to\infty$.
\end{itemize}
In higher dimensions $\R^d, d\ge 3$ one can consider more general kernels
\be
K(x-y)=\left\{
\begin{array}{lll}
\log\frac1{|x-y|}, &d=2, \\
\frac1{|x-y|^{d-2}}, &d\ge 3,
\end{array}
\right.
\ee
with  $\Gamma$ being a Lyapunov surface in $\R^d$ and define the energy 
as follows
\be
I[\mu]=\int\int K({x-y})d\mu(x)d\mu(y)+\int Qd\mu.
\ee
In this note we mostly confine ourselves with quadratic  potentials $Q(x)=|x|^2$ in $\R^2$, although all our results remain valid for more general $Q$ satisfying $\bf (H1)-(H3)$.
Furthermore, our main result on global $L^2$ estimate of the gradient of the equilibrium potential with kernel $K(x-y)=|x-y|^{-d}$ remains valid in in $\R^d, d\ge3$, see Theorem \ref{thm-2}.
\smallskip 

The functional $I[\mu]$, with $Q=|x|^2, d=2$, arises in the description of the 
convergence of the spectral measure of square $N\times N$ matrices with complex independent, standard Gaussian entries (i.e., the Ginibre ensemble) as $N\to \infty$.
In case when there are no constraints imposed on the eigenvalues, it is well known that the 
eigenvalues spread evenly  in the ball of radius $\sqrt N$, and after renormalization by a factor $\frac1{\sqrt N}$ the normalized spectral measure converges to the characteristic function of the unit disc. This is known as the circular law \cite{Ginibre}, \cite{ASZ}. In this context the functional 
$I$ is used to prove large deviation principles for the spectral measure.

\smallskip 

If one demands that the eigenvalues are real (i.e. when $a=1, \Gamma=\R$) we get the so called semicircle law.
More generally, one can demand that a portion of  eigenvalues is contained in a prescribed set $\Gamma$. 
This is considered in \cite{ASZ} when a portion of eigenvalues are contained in an open bounded subset 
 of $\R^2$ and in \cite{Ginibre} when $\Gamma$ is a line. 
 These problems can be related to the thin obstacle and obstacle problems respectively.
 The key step in proving this is to establish  $H^1_{loc}(\R^2)$ estimates for the logarithmic potential 
 $${\ua}=K*\mu_a$$ of the corresponding equilibrium measure.
 The aim of this note is to show that the thin obstacle problem is  well-posed in $H^1(\R^2)$
by showing that in fact $U^{\mu_a}\in H^1(\R^2)$, see Theorem \ref{thm-2}. This improves the previous results in \cite{ASZ} and \cite{Ginibre}.

\smallskip 
The paper is organized as follows: In the next section we prove the  existence and uniqueness of the equilibrium measure $\mu_a$ minimizing the energy $I[\mu]$.
In section 3 we discuss some  basic properties of $\mu_a$. In particular we show that 
there are two positive constants  $A_\Gamma$ and $A_0$ such that 
$2\ua+Q=A_\Gamma$ on $\supp \mu_a\cap \Gamma$ and $2\ua+Q=A_0$ on 
$\supp\mu_a\setminus \Gamma$. Furthermore,  $A_\Gamma>A_0$.
This fact will be used later to show that $\supp\mu_a\setminus \Gamma$ and $\supp\mu_a\cap \Gamma$ are disjoint.

Our main result Theorem \ref{thm-2} is contained in section 4. To prove it we study the   Fourier transformations of $\ua$ and $\mu_a$. It leads to some integral identity involving  Bessel functions. This approach is based on a method of L. Carleson \cite{Carleson}.
Finally, combining the results obtained, in section 5 we show that $\ua$ solves the obstacle problem where the obstacle is given by 
\be
\psi(x)=\left\{
\begin{array}{lll}
\frac12(A_\Gamma-|x|^2) &\hbox{if}\ x\in \Gamma,\\
\frac12(A_0-|x|^2) &\hbox{if}\ x\in \R^2\setminus \Gamma.
\end{array}
\right.
\ee

\section{Existence of minimizers}
In this section we show the existence of a unique equilibrium measure. 
\begin{theorem}\label{thm-1}
Suppose $d=2, \Gamma\subset \R^2$ is a regular $C^{1, \alpha}$ smooth planar curve without self-intersections. 
There is a unique minimizer $\mu_a\in \M_a$ of $I[\mu]$ such that 
\[I[\mu_a]=\inf_{\mu\in \M_a}I[\mu].\]
\end{theorem}
\begin{proof}
Observe that the uniqueness follows from the convexity of 
$\M_a$ and can be proved as in \cite{Ginibre}. Moreover, $I[\mu]$ is also 
semicontinuous. Thus, we have to show that $I[\mu]$ is bounded by below for all $\mu\in \M_a$
and there is at least one $\mu_0$ such that $I[\mu]$ is finite. 
The lower bound follows as in the proof of Theorem 1.3 (a) p. 27 \cite{Saff}.

It remains show that the $\inf_{\mu\in \M_a} I[\mu] <\infty$. Let $\chi_D$ denote the characteristic function of the set $D$ and take
\be\nonumber
\mu=a\frac1L\H^1\with(\Gamma\cap\Om)+(1-a)\frac1{|B|}\chi_B 
\ee
where $B=B_\rho(z)=\{x\in \R^2 \ :\ |x-z|<\rho\}$ with small $\rho$ such that $B\subset \Om$, $\Om\subset \R^2$ is a compact, $L=\H^1(\Gamma\cap\Om)>0$,   and
$\dist(\Gamma, B)>0$. Observe that for this choice of $\mu$ we have 
\be\nonumber
\int_\Om\log\frac1{|x-y|}d\mu(x)=\frac1L\int_0^L\log\frac1{|\Gamma(t)-y|}dt+
\frac1{|B|}\int_B\log\frac1{|x-y|}d\mu(x).
\ee
Assuming that $\Gamma$ is given by arc-length parametrization we have for the logarithmic energy
\be\label{hav-1}
\mathcal L[\mu]=\frac{a^2}{L^2}\int_0^L\int_0^L\log\frac1{|\Gamma(t)-\Gamma(s)|}dtds+\frac{2a(1-a)}{L|B|}\int_0^L\int_{B}\log \frac1{|\Gamma(t)-y|}dtdy+\frac{(1-a)^2}{|B|^2}\int_B\int_B\log\frac1{|x-y|}dxdy.
\ee
Since $\dist(\Gamma, B)>0$ then the second integral is bounded.
As for the last integral then after change of variables $x-y=\xi$ we have 
\be\nonumber
\int_{B\rho(z)}\log\frac1{|x-y|}dx=\int_{B_\rho(z-y)}\log\frac1{|\xi|}d\xi\le \int_{B_{2\rho}(0)}\log\frac1{|\xi|}dx<\infty
\ee
where we used $|z-y|\le \rho$ and the fact that $\rho$ is small by construction.

It remains to check that the first integral is finite. Let us fix $s\in [0, L]$
Then we have that 
\begin{eqnarray*}
\int_0^L\log \frac1{|\Gamma(t)-\Gamma(s)|}dt&=&\int_{-s}^{L-s}\log \frac1{|\Gamma(\tau+s)-\Gamma(s)|}d\tau=\\\
&=&\tau\log \frac1{|\Gamma(\tau+s)-\Gamma(s)|}\bigg|_{-s}^{L-s}-\int_{-s}^{L-s} \tau \frac{\dot{\Gamma}(\tau+s)\cdot(\Gamma(\tau+s)-\Gamma(s))}{|\Gamma(\tau+s)-\Gamma(s)|^2}d\tau=\\\nonumber
&=&(L-s)\log\frac1{|\Gamma(L)-\Gamma(s)|}+s\log\frac1{| \Gamma(0)-\Gamma(s)|}-I_0
\end{eqnarray*}
where $I_0$ is the last integral. Using the crude estimate
\begin{eqnarray}\label{hav-2}
|I_0|&\le&\int_{-s}^{L-s}|\tau|\frac{|\dot\Gamma(\tau+s)|}{|\Gamma(\tau+s)-\Gamma(s)|}d\tau=\int_{-s}^{L-s}\frac{|\tau|}{|\Gamma(\tau+s)-\Gamma(s)|}d\tau=\\\nonumber
&=&\int_{[-s, L-s]\setminus(-\delta, \delta)}\frac{|\tau|}{|\Gamma(\tau+s)-\Gamma(s)|}d\tau+\int_{-\delta}^{\delta}
\frac{|\tau|}{|\Gamma(\tau+s)-\Gamma(s)|}d\tau\\\nonumber
&\le& \frac{4L^2}{C_\delta} +\int_{-\delta}^{\delta}
\frac{|\tau|}{|\Gamma(\tau+s)-\Gamma(s)|}d\tau
\end{eqnarray}
because $|\Gamma(\tau+s)-\Gamma(s)|\ge C_\delta$ if $|\tau|\ge \delta$.
 Finally, from $C^{1, \alpha}$ regularity of $\Gamma$ we get 
 \begin{eqnarray}\label{hav-3}
|\Gamma(\tau+s)-\Gamma(s)|&=&|\tau|\left|\int_0^1\dot \Gamma(\sigma\tau+s)d\sigma\right|\ge \\\nonumber
&\ge& |\tau|\left(|\dot\Gamma(s)|-\int_0^1|\dot \Gamma(\sigma\tau+s)-\dot\Gamma(s)|d\sigma\right)\\\nonumber
&\ge&|\tau|\left(1-\delta^\alpha \right).
 \end{eqnarray}
Combining \eqref{hav-3} with \eqref{hav-2} we get 
\be\nonumber
|I_0|\le \frac{4L^2}{C_\delta} +2\delta\left(1-\delta^\alpha \right)<\infty.
\ee
Returning to the first integral  in \eqref{hav-1} we infer 
\begin{eqnarray*}
\int_0^L\int_0^L\log\frac1{|\Gamma(t)-\Gamma(s)|}dtds&\le&\int_0^L \left\{(L-s)\log\frac1{|\Gamma(L)-\Gamma(s)|}+s\log\frac1{| \Gamma(0)-\Gamma(s)|}+\frac{4L^2}{C_\delta} +2\delta\left(1-\delta^\alpha \right)\right\}ds\\\nonumber
&\le&L\left[\frac{4L^2}{C_\delta} +2\delta\left(1-\delta^\alpha \right)\right]+L\log\frac1{C_\delta}+\\\nonumber
&&+\int_{\delta}^{L-\delta}\left\{(L-s)\log\frac1{|\Gamma(L)-\Gamma(s)|}+s\log\frac1{| \Gamma(0)-\Gamma(s)|}\right\}ds\\\nonumber
&\le& C(\delta, L)
\end{eqnarray*}
if we choose $\delta>0$ suitably small.
This finishes the proof for $d=2$.
\end{proof}

\begin{remark}
If $d\ge 3, Q(x)=|x|^2$ then clearly $I[\mu]\ge 0$. The upper estimate for $I[\mu]$ follows from a similar argument if we assume that $\Gamma$ is a Lyapunov surface and take $\mu=a\frac1L\H^{d-1}\with(\Gamma\cap\Om)+(1-a)\frac1{|B|}\chi_B 
$ with $L=\H^{d-1}(\Gamma\cap \Om)$ and $\dist(B, \Gamma)>0$.  Therefore Theorem \ref{thm-1} remains valid for $d\ge 3.$
\end{remark}
\section{Basic properties of minimizers}
In this section we prove some basic properties of the equilibrium measure.
The arguments are along the line of those in \cite{ASZ}.
Therefore, we mostly focus on those aspects of the proofs which are new or differ essentially. 
The results to follow are valid in $\R^d, d\ge 2$ unless otherwise stated.
\begin{lem}\label{lem-equit}
 Let $\mu_a$ be as in Theorem \ref{thm-1}. Then $\mu_a(\Gamma)=a$.
\end{lem}
\begin{proof}
If the claim fails then $\mu_a(\Gamma)>a$. 
Fix $\delta\in(0, a)$ and let $\mu_{a-\delta}$ be the minimizer of $I[\cdot]$ over $\M_{a-\delta}\supset \M_a$.
Form $\mu=(1-\e)\mu_a+\e\mu_{a-\delta}, \e\in[0,1]$.
Clearly, $\mu\in \M_{a}$ if we choose $\e\delta$ sufficiently small because 
\[\mu(\Gamma)>a+[\mu_a(\Gamma)-a]-\e\delta.\]
Consequently, we have from the strict convexity of $I$
\begin{eqnarray*}
I[(1-\e)\mu_a+\e\mu_{a-\delta}]&<&(1-\e)I[\mu_a]+\e I[\mu_{a-\delta}]=I[\mu_a]+\e\left(I[\mu_{a-\delta}]-I[\mu_a]\right)\\\nonumber
&\le& I[\mu_a]
\end{eqnarray*}
which is in contradiction with the fact that $\mu_a$ is a minimizer.
\end{proof}

Observe that the Fr\'echet derivative of $I[\mu]$ is $2\ua+Q$
where 
\[\ua(y)=\int K(x-y)d\mu_a(x).\] 

It is convenient to 
consider variations of the equilibrium measure in terms of affine combinations.
More precisely, let $\mu_\e=(1-\e)\mu_a+\e\nu, \nu\in \M_a, \e\in [0,1]$,  then 
by direct computation we have that 
\begin{eqnarray}\label{blya-1}
I[\mu_\e]
&=&(1-\e)^2\int\int K(x-y)d\mu_a(x)d\mu_a(y)\\\nonumber
&&+2\e(1-\e)\int\int K(x-y)d\mu_a(x)d\nu(y)+\e^2\int \int K(x-y)d\nu(x)d\nu(y)\\\nonumber
&&+(1-\e)\int Qd\mu_a+\e\int Qd\nu\\\nonumber
&=&I[\mu_a]+\e\left( 
2\int\int K(x-y)d\mu_a(x)d(\nu(y)-\mu_a)+\int Qd(\nu-\mu_a)\right)+O(\e^2)=\\\nonumber
&=&I[\mu_a]+\e\int(2U^{\mu_a}+Q)d(\nu-\mu_a)+O(\e^2).
\end{eqnarray}
Since $\mu_a$ is the minimizer then $I[\mu_a]\le I[\mu]$, and after sending $\e\to 0$ it follows 
that 
\begin{equation}
\int(2U^{\mu_a}+Q)d(\nu-\mu_a)\ge 0, \quad \forall \nu\in \M_a.
\end{equation}

 \begin{lem}\label{lem-hor}
 Let $A_\Gamma=\frac1{a}\int_{\Gamma}(2U^{\mu_a}+Q)d\mu_a$ then
 quasi everywhere
 
 \begin{eqnarray}\label{A-gam}
 2U^{\mu_a}+Q
\begin{array}{lll}
=A_\Gamma &\hbox{on}\ \ \Gamma\cap \supp\mu_a, \\
\ge A_\Gamma & \hbox{on} \ \ \Gamma.
\end{array}
 \end{eqnarray}
 Similarly, let us denote $A_0=\frac1{1-a}\int_{\R^2\setminus\Gamma} (2U^{\mu_a}+Q)d\mu_a$
 then 
  \begin{eqnarray}\label{A-zero}
 2U^{\mu_a}+Q
\begin{array}{lll}
=A_0 &\hbox{on}\ \ \supp\mu_a\setminus \Gamma,\\
\ge A_0 & \hbox{on} \ \  \R^2\setminus (\supp\mu_a\setminus \Gamma). 
\end{array}
 \end{eqnarray}
 Furthermore, 
 \begin{equation}\label{const-ineq}
 A_\Gamma>A_0.
 \end{equation}

 \end{lem}
\begin{proof}
We first prove \eqref{A-gam}. Suppose that there is a set capacitable $E$ of positive capacity such that 
$\Gamma\cap E$ has zero capacity and 
\[2U^{\mu_a}+Q<A_\Gamma-\delta\ \  \hbox{q.e. on} \ \  E \]
for some positive $\delta.$
Let $\mu_E$ be the equilibrium measure of $E$ and form $\nu=\mu_a\with (\R^2\setminus\Gamma)+a\mu_E$.
Clearly $\nu\in \M_a$. Therefore, in view of \eqref{blya-1} for the measure $\mu_\e=\e \mu_a+(1-\e)\nu\in \M_a$ we get 

\begin{eqnarray}
I[\mu_\e]&=& I[\mu_a]+\e\left( 
2\int\int K(x-y)d\mu_a(x)d(\nu(y)-\mu_a)+\int Qd(\nu-\mu_a)\right)+O(\e^2)\\\nonumber
&=&I[\mu_a]+\e 
\int_{\Gamma}(2U^{\mu_a}+Q)d(a\mu_E-\mu_a)+O(\e^2)\\\nonumber
&=&I[\mu_a]+\e\left(a\int_{\Gamma}(2U^{\mu_a}+Q)d\mu_E-a A_\Gamma\right)+O(\e^2)\\\nonumber
&<&I[\mu_a]-a \e\delta+O(\e^2)\\\nonumber 
&<&I[\mu_a]
\end{eqnarray}
if $\e$ and $\delta$ are sufficiently small. This will be in contradiction with the fact that 
$\mu_a$ is the minimizer.
Thus we have proved that $2U^{\mu_a}+Q\ge A_\Gamma$ q.e. on $ \Gamma$.

Next we show that 
on $\supp\mu_a\cap \Gamma$ we have $2U^{\mu_a}+Q=A_\Gamma$ q.e. Indeed, from the  definition of 
$A_\Gamma$ it follows
\[aA_\Gamma=\int_{\Gamma}(2U^{\mu_a}+Q)d\mu_a\ge aA_\Gamma\]
where the last inequality follows from the first inequality in \eqref{A-gam}.
The proof of \eqref{A-zero} is similar.
In order to prove the last claim $A_\Gamma>A_0$ we first observe that there exists a 
measure $\nu\in\M_a$ such that 
\begin{itemize}
\item $a>\nu(\Gamma)$,
\item $I[\nu]\le I[\mu_a]$.
\end{itemize}
First notice that $\M_a\subset \M_{a-\delta}$ for $\delta\in (0, a)$.
Fix such $\delta>0$ and let $\mu_{a-\delta}$ be the minimizer of $I[\cdot]$ over $\M_{a-\delta}$.
Then by Lemma \ref{lem-equit} $\mu_{a-\delta}(\Gamma)=a-\delta<a$
and $I[\mu_{a-\delta}]=\inf_{\M_{a-\delta}}I[\mu]\le I[\mu_a]=\inf_{\M_a}I[\mu]$. Therefore one can take $\nu=\mu_{a-\delta}$.

From the strict convexity of $I$ it follows that 
\[I[\nu]>I[\mu_a]+\langle DI[\mu_a], \nu-\mu_a\rangle\]
where $DI[\mu]=2U^{\mu}+Q$ is the Fr\'echet derivative of $I[\mu]$.
Therefore, from the properties of $\nu$ we infer 
\be\label{blya-3}
0\ge I[\nu]-I[\mu_a]>\langle DI[\mu_a], \nu-\mu_a\rangle
\ee
or equivalently 
\[\langle 2U^{\mu_a}+Q, \nu-\mu_a\rangle<0.\]
On the other hand 
\be\label{blya-2}
\int(2U^{\mu_a}+Q)d\mu_a=a A_\Gamma+(1-a)A_0\ee
while
\[\int(2U^{\mu_a}+Q)d\nu=\int_\Gamma(2U^{\mu_a}+Q)d\nu+\int_{\R^2\setminus
\Gamma}(2U^{\mu_a}+Q)d\nu\ge \nu(\Gamma)A_\Gamma+\nu(\R^2\setminus\Gamma)A_0.\]
This together with \eqref{blya-2}, \eqref{blya-3}
 yields
\[a A_\Gamma+(1-a)A_0>\nu(\Gamma)A_\Gamma+(1-\nu(\Gamma))A_0\Rightarrow A_0(\nu(\Gamma)-a)>A_\Gamma(\nu(\Gamma)-a).\]
Finally, the property  $\nu(\Gamma)<a$ implies that $A_\Gamma>A_0$.
\end{proof}

\begin{cor}\label{cor-comp-supp}
$\supp\mu_a$ is compact.
\end{cor}
\begin{proof}
If $d\ge 3$ then $K(x-y)\ge0$, hence by Lemma \ref{lem-hor} for $x\in \supp\mu_a$ we have 
\be
\max(A_\Gamma, A_0)\ge 2U^{\mu_a}(x)+Q(x)\ge Q(x)\to \infty \ \ \ \  \hbox{if}\ \ \ |x|\to \infty
\ee
which is a contradiction.
If $d=2$ then from the triangle inequality we get that 
\be
K(x-y)\ge -\log|x|-\log\left(1+\frac{|y|}{|x|}\right).
\ee
Consequently, for $x\in \supp\mu_a$ 
\begin{eqnarray*}
\max(A_\Gamma, A_0)&\ge& 2U^{\mu_a}(x)+Q(x)\ge Q(x)-2\log|x|-\int\log\left(1+\frac{|y|}{|x|}\right)d\mu_a\\\nonumber
&=&Q(x)-2\log|x|+O(1)\to \infty \ \ \ \  \hbox{if}\ \ \ |x|\to \infty
\end{eqnarray*}
for sufficiently large $|x|$, where the last inequality follows from \eqref{Carleson} and $\int Qd\mu_a<I[\mu_a]<\infty$. Since $Q=|x|^2$ (of for the general case from the hypotheses on $Q$ $\bf(H1)-(H3)$) it again follows that 
$\supp \mu_a$ is bounded.
\end{proof}

\section{Global $L^2$ estimates for $U^{\mu_a}$ and $\na U^{\mu_a}$}
Our main result is contained in the following 
\begin{theorem}\label{thm-2}
Let $U^{\mu_a}(y)=\int K(x-y)d\mu_a$,  if $d\ge 3$ then $\na \ua\in L^2(\R^d)$. 
If $d=2$ then $U^{\mu_a}\in H^1(\R^2)$.
Furthermore, there holds
\be 
 \|\ua\|_{ H^1(\R^2)}\le C\mathcal E[\mu_a].
 \ee
 Here $\mathcal E[\mu]$ is the energy of $\mu$ defined as $\int\int K(x-y)d\mu(x)d\mu(y)$. 
\end{theorem}
\begin{remark}
It is shown in \cite{Carleson} that $\mathcal E[\mu]>0$ for any probability measure $\mu$
and $d\ge2$.
In fact, this can be seen from the proof to follow (see also Corollary \ref{cor-posit}).
\end{remark}

\begin{proof}
The case $d\ge 3$ follows from Lemma 1.6 p. 92 \cite{Landkof} (see also Lemma 17 p. 95), which assert that 
\[\frac{\p \ua(x)}{\p x_i}=\int \frac{\p K(x-y)}{\p x_i}d\mu_a\]
almost everywhere  and morover
\[\frac1{4\pi^2}\int_{\R^d}|\na \ua|^2\le \int\int K(x-y)d\mu_a(x)d\mu_a(y)=\mathcal E[\mu_a].\]
The case of the logarithmic potential follows from a modification of the argument by L. Carleson 
\cite{Carleson} Lemma 3 page 22. We begin with computing the Fourier transformation of $K$.
Note that since $\supp\mu_a$ is compact  we can assume that $K(r)=0$ for $r\ge r_0$
for some fixed $r_0>0$. 
We have 
\begin{eqnarray*}
\F K(\xi)&=&\int K(x)e^{-2\pi i\langle x, \xi\rangle}dx=\int K(x)e^{-2\pi i\langle x|\xi|, \frac\xi{|\xi|}\rangle}dx\\\nonumber
&=&\frac1{4\pi^2|\xi|^2}\int K\left(\frac y{2\pi |\xi|}\right)e^{i\langle y, \frac\xi{|\xi|}\rangle}dx.
\end{eqnarray*}
Let us denote $K_0(y)=K\left(\frac y{2\pi |\xi|}\right)$ and define 
\[F(\eta)=\int K_0(y)e^{i\pi\langle y, \eta\rangle}, \quad \eta=\frac\xi{|\xi|}.\]
From Lemma 2 p. 21 \cite{Carleson} it follows that there is a universal constant $c_1$ such that 
\[F(\eta)=c_1\int_0^\infty K_0(r)J(r)rdr, \quad |\eta|=1\]
where $J$ is the Bessel function 
\be
J(r)=-J''(r)-\frac{J'(r)}r, \quad J(0)=1, J'(0)=0, \quad J(r)<1, r\not =0.
\ee
Therefore $F(\eta)$ can be further simplified as follows 
\begin{eqnarray} 
F(\eta)&=&-c_1\int_0^\infty K_0(r)(rJ(r))'dr=\\\nonumber
&=&c_1\int_0^{2\pi |\xi|r_0} rJ'(r)K_0'(r)dr
\end{eqnarray}
because from the definition of $K_0$ we have $\supp K_0\subset [0, 2\pi|\xi|r_0]$.
Moreover, $K_0'(r)=-\frac1r$ hence 
\be
F(\eta)=c_1(1-J(2\pi |\xi|r_0)).
\ee
Consequently, 
\be 
\F K(\xi)=\frac{c_1}{4\pi^2|\xi|^2}(1-J(2\pi |\xi|r_0)).
\ee
Next we restrict $\mu_1=\mu_a\with \mathcal C$ where $\mathcal C\subset \supp\mu_a$
is a compact such that $U^{\mu_1}$ is continuous. Observe that 
$\int \ua d\mu_a$ is finite hence $\ua$ is finite $\mu_a$ almost everywhere. 
By Theorem 1.8  p. 70 \cite{Landkof} for every $\e>0$ small there is a  restriction of 
$\mu_a$ such that 
\[0\le \int\mu_a-\int\mu_1<\e.\]
Note that if $\tau =\mu_a-\mu_1$ then we have 
\[|\mathcal E[\mu_a]-\mathcal E[\mu_1]|=\left|\int U^{\mu_a-\mu_1}d\mu_a+\int U^{\mu_a-\mu_1}d\mu_1\right|=\left|\int (\ua+U^{\mu_1})d\tau\right|=O(\e).\]
Let $\phi_n(y)=n^{\frac d2}e^{-n\pi |y|^2}$ be the sequence of normalised 
Gaussian kernels. It is well-known that $\phi_n$ is a mollification kernel for every $n\in \mathbb N$ and moreover $\F{\phi_n}=e^{-\frac{\phi|\xi|^2}n}$. From the Parseval relation 
\be
\int(\phi_n*U^{\mu_1})d\mu_1=\int\F{\phi_n}\F K|\F{\mu_1}|^2.
\ee 
If we first send $n\to \infty$ and then $\e\to 0$ to conclude the identity
\be
\mathcal E[\mu_a]=\int\F K |\F{\mu_a}|^2.
\ee
 On the other hand $\F\ua=\F K\F{\mu_a}$,  which yields
\begin{eqnarray}\label{blya-4}
\mathcal E[\mu_a]&=&\int\F K (\xi)\frac{|\F{U^{\mu_a}}(\xi)|^2}{|\F K(\xi)|^2}d\xi\\\nonumber
&=&\int\frac{4\pi^2|\xi|^2}{c_1(1-J(2\pi r_0|\xi|))}|\F{U^{\mu_a}}(\xi)|^2d\xi\\\nonumber
&=&\int_{|\xi|<\delta}+\int_{|\xi|\ge \delta}.
\end{eqnarray}
Using the expansion $J(t)=\sum_{s=0}^\infty \frac{(-1)^s}{(s!)^2}\left(\frac t2\right)^{2s}=1-\frac{t^2}4+\frac{t^4}{64}+\dots$
we see that \[\frac{4\pi^2|\xi|^2}{c_1(1-J(2\pi r_0|\xi|))}=\frac{1}{r_0^2c_1}
\frac4{(1-\frac{(2\pi r_0|\xi|)^2}{16}+\dots)}\]
hence the first integral is bounded below by $C(\delta)\frac{1}{r_0^2c_1}\int_{|\xi|<\delta}|\F{U^{\mu_a}}(\xi)|^2d\xi$ for sufficiently small $\delta>0$.
As for the second integral, we have 
\be 
\int_{|\xi|\ge\delta}\frac{4\pi^2|\xi|^2}{c_1(1-J(2\pi r_0|\xi|))}|\F{U^{\mu_a}}(\xi)|^2d\xi\ge 
\frac{4\pi^2\delta^2}{c_1}\int_{|\xi|\ge\delta}|\F{U^{\mu_a}}(\xi)|^2d\xi.
\ee
Combining we see that $\F\ua\in L^2(\R^2)$ which, after we apply Parseval's relation again, yields
 $\ua \in L^2(\R^2)$ and
 \be 
 \|\ua\|_{ L^2(\R^2)}\le C\mathcal E[\mu_a].
 \ee
 To finish the proof we use that $4\pi^2|\xi|^2|\F \ua|^2=|\F{\na\ua}|^2$ which together with 
 \eqref{blya-4} implies that 
 \be
 \mathcal E[\mu_a]
=\int\frac{1}{c_1(1-J(2\pi r_0|\xi|))}|\F{\na U^{\mu_a}}(\xi)|^2d\xi\ge 
\frac1{c_1}\int|\F{\na U^{\mu_a}}(\xi)|^2d\xi\
 \ee
 which finishes the proof.
\end{proof}

\begin{cor}\label{cor-posit}
Let $\mu_a$ be as in Theorem \ref{thm-1}. Then there holds 
\be\label{Carleson}
\mathcal E [\mu_a]=\int U^{\mu_a}d\mu_a>0.
\ee
\end{cor}

\section{The  thin obstacle problem}

From the $H^1(\R^2)$ estimate for $\ua$ it follows that $\ua$ is a solution  to some  
variational inequality, and hence $\ua$ can be interpreted as a solution to an obstacle problem
with a combination of both thin (on $\Gamma$) and "thick" obstacles (on $\R^2\setminus\Gamma$). It is convenient to define the obstacle as follows
\be
\psi(x)=\left\{
\begin{array}{lll}
\frac12(A_\Gamma-|x|^2) &\hbox{if}\ x\in \Gamma,\\
\frac12(A_0-|x|^2) &\hbox{if}\ x\in \R^2\setminus \Gamma.
\end{array}
\right.
\ee

\begin{lem}
Let $\ua$ be the logarithmic potential of $\mu_a$ and define 
$$\mathcal K=\left\{v\in H^1_{loc}(\R^2)\ s.t. \ v-\ua \ \hbox{has bounded support in}\ \R^2,\ v\ge \psi\right\}.$$
Then $\ua$ solves the following obstacle problem:
\[\int\na \ua \na(v-\ua)\ge0, \quad \forall v\in \mathcal K.\]
\end{lem}
The proof is the same as in \cite{ASZ}.

\begin{cor}\label{cor-sep}
$\dist(\Gamma, \supp(\mu_a\setminus\Gamma))>0$.
\end{cor}
\begin{proof}
This follows from the estimate $A_\Gamma>A_0$. Indeed, let us assume that 
$x_0\in \Gamma\cap \supp \mu_a$ and there is a sequence $\{x_k\}_{k=1}^\infty, x_k\in \supp\mu_a\setminus \Gamma$ such that $\lim_{k\to \infty}x_k\to x_0$. Using the lower semicontinuity of $\ua$ (see Lemma 1 p.15 \cite{Carleson}) we see that 
\be\label{mek}
\frac12(A_0-|x_0|^2)=\liminf_{x_k\to x_0}\ua(x_k)\ge \ua(x_0). 
\ee
Let $\rho>0$ be such that $\{x_k\}\subset B_{\rho}(x_0)$. If $\rho$ is small then 
$\Gamma $ divides $B_\rho(x_0)$ into two parts $D^+$ and $D^-$. To fix the ideas let us suppose that $D^+$ contains a subsequence $\{x_k\}$. Let $h$ be the harmonic function in $D^+$ such that 
$h=\psi$ on $\p D^+$.
Observe that $h$ is continuous at $x_0$ because $\Gamma\in C^{1, \alpha}$. Since $\ua$ is   superharmonic and on 
$\p D^+$ we have $\ua\ge \psi=h$ then the comparison principle implies that 
\be\label{erku}
\ua(x_0)\ge h(x_0)=\frac12(A_\Gamma-|x_0|^2).
\ee
Combining \eqref{mek} and \eqref{erku} we see that $A_0\ge A_\Gamma$
 which is a contradiction in view of \eqref{const-ineq}.
 \end{proof}

From Corollary \ref{cor-sep} it follows that near $\Gamma$ the potential $\ua$ is a solution to a thin obstacle 
problem in the following sense, see \cite{Friedman} p. 108: 
\begin{equation}
\left.
\begin{array}{lll}
\ua\ge \frac12(A_\Gamma-Q)\\
\frac{\p \ua}{\p n^+}+ \frac{\p \ua}{\p n^-}\ge 0\\
\left(u-\frac12(A_\Gamma-Q)\right)\left(\frac{\p \ua}{\p n^+}+ \frac{\p \ua}{\p n^-}\right)=0
\end{array}
\right\}\hbox{on}\ \Gamma
\end{equation}
where $n^\pm$ are the outward normals on the $\Gamma$ corresponding to the domains that $\Gamma$ separates. In particular, if $\Gamma$ is $C^3$ regular then 
$\ua$ is $C^{1, \alpha}$ up to $\Gamma$ from each of its side, see Theorem 11.4 p.111 \cite{Friedman}.

A particular case is $\Gamma=\R$ \cite{Ginibre}. Using a simple symmetrization argument (see e.g. \cite{Hayman} p. 119 Theorem 4.6) we can show that the  
potential $\ua$ is symmetric w.r.t. the real line and hence we get the Signorini problem  near $\R$ \cite{Friedman} p. 111. 

One can make the connections with the obstacle problem more explicit by using the 
$H^1(\R^2)$ estimate in Theorem \ref{thm-2} and transforming the 
energy $I[\mu_a]$. Let $R>0$ be fixed then using the divergence theorem 
\begin{eqnarray}\label{ereq}
\int_{B_R}\ua d\mu_a&=&-\frac1{2\pi} \int_{B_R}\ua\Delta \ua =\\\nonumber
&=&\frac1{2\pi} \int_{B_R}|\na \ua |^2-\frac1{2\pi}\int_{\p B_R}\ua \p_n \ua.
\end{eqnarray}
For a.e. $R>0$ the last integral can be estimated as follows 
\[
\left|\int_{\p B_R}\ua \p_n \ua\right|\le \int_{\p B_R}|\ua| |\na  \ua |\le \int_{\p B_R}|\ua|^2+ |\na \ua
|^2.\]
From Theorem \ref{thm-2} and Fubini's theorem it follows that 
\[
\int_{R^2}(|\ua|^2+|\na \ua|^2)=\int_0^\infty\int_{\p B_R}(|\ua|^2+ |\na \ua
|^2)dR.
\]
Consequently,  
\[\int_{\p B_R}|\ua|^2+ |\na \ua
|^2\to 0\quad  R\to \infty\]
and we infer from \eqref{ereq} that 
\[\int_{\R^2}\ua d\mu_a=\frac1{2\pi} \int_{\R^2}|\na \ua |^2.\]
Recalling that by Corollary \ref{cor-comp-supp} $\supp \mu_a\subset B_{r_0}$ for some $r_0>0$ 
and using the divergence theorem again we conclude 
\begin{eqnarray}
\int_{B_{r_0}}|x|^2d\mu_a&=&-\frac1{2\pi}\int_{B_{r_0}}|x|^2\Delta \ua =-\frac1{2\pi}\int_{B_{r_0}}
\ua \Delta |x|^2+\frac1{2\pi}\int_{\p B_{r_0}}(2r_0 \ua -r_0^2\p_n\ua)\\\nonumber
&=&-\frac2{\pi }\int_{B_{r_0}}\ua +\frac{r_0}{\pi}\int_{\p B_{r_0}}\ua +r_0^2.
\end{eqnarray}
Combining these we have that the energy can be rewritten in terms of $\ua$ in the following form
\[I[\mu_a]=\frac1{2\pi} \int_{\R^2}|\na \ua |^2-\frac2{\pi }\int_{B_{r_0}}\ua +\frac{r_0}{\pi}\int_{\p B_{r_0}}\ua +r_0^2.
\]


\end{document}